\documentclass{gtart_a}
\pdfoutput=1

\title{Global rigidity for totally nonsymplectic Anosov $\mathbb{Z}^k$ actions} 

\author{Boris Kalinin}
\givenname{Boris}
\surname{Kalinin}
\address{Department of Mathematics and Statistics\\
 University of South  Alabama\\\newline
Mobile, AL 36688 \\USA}
\email{kalinin@jaguar1.usouthal.edu}
\urladdr{}

\author{Victoria Sadovskaya}
\givenname{Victoria}
\surname{Sadovskaya}
\email{sadovska@jaguar1.usouthal.edu}
\urladdr{}

\volumenumber{10}
\issuenumber{}
\publicationyear{2006}
\papernumber{24}
\lognumber{0653}
\startpage{929}
\endpage{954}

\doi{}
\MR{}
\Zbl{}

\arxivreference{math/0602175}  
\arxivpassword{}   


\received{8 September 2005}
\revised{11 June 2006}
\accepted{5 June 2006}
\published{24 July 2006}
\publishedonline{24 July 2006}
\proposed{Benson Farb}
\seconded{David Gabai, Leonid Polterovich}
\corresponding{}
\editor{}
\version{}

\AtBeginDocument{\let\bar\wbar\let\tilde\wtilde}
\AtBeginDocument{\def \a{\alpha}\def \L{\mathcal{L}}\def \H{\mathcal{H}}\def \b{\beta}}

\def\co{\mskip0.5mu\colon\thinspace}

\makeatletter
\def\cnewtheorem#1[#2]#3{\newtheorem{#1}{#3}[section]
\expandafter\let\csname c@#1\endcsname\c@theorem}

\theoremstyle{plain}
      \newtheorem{theorem}{Theorem}[section]
      \cnewtheorem{lemma}[theorem]{Lemma}
      \cnewtheorem{corollary}[theorem]{Corollary}
      \cnewtheorem{proposition}[theorem]{Proposition}
\theoremstyle{remark}  
      \cnewtheorem{remark}[theorem]{Remark}
\theoremstyle{definition}       
      \cnewtheorem{conjecture}[theorem]{Conjecture}
      \cnewtheorem{definition}[theorem]{Definition}     
      
\makeatother  

\numberwithin{equation}{section}       

\def \Rk{\mathbb R^k}
\def \Rn{\mathbb R^n}
\def \Rm{\mathbb R^m}
\def \Tn{\mathbb T^n}
\def \Zk{\mathbb Z^k}

\def \Es{E^{s}}

\def \Ws{W^{s}}
\def \Ci{C^{\infty}}
\def \bH{\bar H}
\def\Id{\text{Id}}
\def\tE{\tilde E}

\def \ta{\tilde \alpha}

\def \la{{\lambda}}

\def \n{\mathfrak{n}}
\def \E{\mathfrak{e}}
\def \e{\varepsilon}
\def \A{\mathcal{A}}
\def \K{\mathcal{K}}
\def \M{\mathcal{M}}
\def \N{\mathcal{N}}

\makeop{dist}
\makeop{ker}

\begin{document}

\begin{asciiabstract}
We consider a totally nonsymplectic Anosov action of Z^k
which is either uniformly quasiconformal or pinched on each coarse 
Lyapunov distribution. We show that such an action on a torus is 
C^\infty--conjugate to an action by affine automorphisms. We also 
obtain similar global rigidity results for actions on an arbitrary 
compact manifold assuming that the coarse Lyapunov foliations are 
jointly integrable.
\end{asciiabstract}

\begin{abstract} 
We consider a totally nonsymplectic (TNS) Anosov action of $\mathbb Z^k$
which is either uniformly quasiconformal or pinched on each coarse Lyapunov 
distribution. We show that such an action on a torus is $C^{\infty}$--conjugate to an 
action by affine automorphisms.  We also obtain similar global rigidity 
results for actions on an arbitrary compact manifold assuming that the 
coarse Lyapunov foliations are topologically jointly integrable.
\end{abstract}

\maketitle


\section{Introduction}

In this paper we consider Anosov actions of $\Zk$, $k\ge 2$, on a compact smooth
manifold. For certain classes of these actions, we establish global rigidity, ie,
the existence of a smooth conjugacy to an algebraic model. 

The classification of Anosov systems is one of the central problems
in smooth dynamics. A long standing conjecture on \textit{topological\/} classification
is that any Anosov diffeomorphism is topologically conjugate to a hyperbolic 
automorphism of a torus or, more generally, an infranilmanifold. This
was established under the assumption that the underlying manifold is an 
infranilmanifold, or that the diffeomorphism is of codimension one (Franks
\cite{Fr}, Manning \cite{M}, Newhouse \cite{N}). 
A smooth classification of Anosov diffeomorphisms is not feasible. Indeed,
even when two Anosov diffeomorphisms are topologically conjugate, the 
conjugacy is typically only H\"older continuous. 

In contrast to a single Anosov diffeomorphism, it may be possible to classify 
higher rank Anosov actions up to a \textit{smooth\/} conjugacy. 
These are the actions of higher rank groups with at least one element acting 
normally hyperbolically to the orbit foliation. The study of rigidity for these actions
originated from Zimmer's conjecture that the standard action of $SL(n,\Z)$ on $\Tn$, 
$n>2$, is locally rigid, ie, any $C^1$--small perturbation is smoothly conjugate to the 
original action. The smoothness of the conjugacy was established using the action 
of a diagonalizable subgroup isomorphic to $\Z ^{n-1}$ (Hurder \cite{H1},
Katok and Lewis \cite{KL2}). 
This prompted the study of higher rank \textit{abelian\/} Anosov actions.

Reducible abelian actions can be obtained as products of Anosov diffeomorphisms
or flows. Important irreducible examples arise from natural algebraic constructions. 
They include $\Zk$ actions by automorphisms of tori and infranilmanifolds, and 
$\Rk$ actions by left translations on homogeneous spaces  and  biquotients. 
These actions exhibit strong rigidity phenomena such as scarcity of invariant 
measures and cocycle rigidity (see surveys in Lindenstrauss \cite{Lind04} and 
Nitica--T\"{o}r\"{o}k \cite{NT}).

It is conjectured that all irreducible Anosov actions of $\Zk$ and $\Rk$, $k\ge 2$, 
are \textit{smoothly\/} conjugate to algebraic actions. This conjecture is supported by the local 
rigidity results of Katok and Spatzier. They showed that any $C^1$--small perturbation
of a higher-rank algebraic Anosov action with semisimple linear part is smoothly
conjugate to the original action \cite{KS97}.  There are few results, however, toward 
proving the conjecture. So far, all known results are for actions with one-dimensional 
coarse Lyapunov distributions. Such actions are often called Cartan actions. Coarse 
Lyapunov distributions are the finest nontrivial intersections of stable distributions 
of Anosov elements of the action. Katok and Lewis established global rigidity for actions 
of certain maximal sets of commuting diffeomorphisms of a torus \cite{KL2}. 
In this case, the coarse Lyapunov distributions were one-dimensional stable 
distributions of some Anosov elements of the action. 
Recently,  the first author obtained a smooth classification of certain classes of 
continuous and discrete actions of rank $k \ge 3$ on arbitrary
manifolds in the joint work with Spatzier \cite{KaSp}.

We consider Anosov $\Zk$ actions with higher-dimensional  coarse Lyapunov 
distributions under various assumptions on the relation between slow and fast 
expansion/con\-traction rates on these distributions. Our approach is different from 
those of \cite{KaSp} or Katok and Lewis \cite{KL2}, it is based on the study of holonomy maps and affine 
structures on coarse Lyapunov foliations. For actions on tori, in contrast to 
\cite{KL2}, we do not rely directly on the topological conjugacy given by the
topological classification. This allows us to obtain results on an arbitrary manifold,
provided that the coarse Lyapunov foliations are topologically jointly integrable. 
Also, we assume that the action is  totally nonsymplectic (TNS), 
ie, any pair of coarse Lyapunov distributions is contracted by some element of the action. 
We note that this assumption is satisfied by the actions in  \cite{KL2} and by the discrete
actions in \cite{KaSp}.

We state our main results in \fullref{Main Results} and prove them in
Sections \ref{Proof Main QC} and \ref{Proof Main pinching}. In \fullref{Preliminaries} we discuss the  structure of $\Zk$ Anosov actions and
give necessary definitions.

The first author was supported in part by NSF grant DMS-0140513, 
and the second author was supported in part by NSF grant DMS-0401014.


\section{Preliminaries}   \label{Preliminaries}


\subsection{Lyapunov exponents for $\Zk$ actions}

In this section we describe the Multiplicative Ergodic Theorem and related
notions for $\Zk$ actions (see Kalinin and Katok \cite{KaK} for more details).
Let  $f$ be a diffeomorphism of a compact manifold $\M$ preserving an 
ergodic probability measure $\mu$. By Oseledec Multiplicative Ergodic 
Theorem, there exist finitely many numbers $\chi _i $ and a measurable 
splitting of the tangent bundle $T\M = \bigoplus$$E _i$ on a set of full measure 
such that the forward and backward Lyapunov exponents of $v \in  E_i $ 
are $\chi _i $.  This splitting is called Lyapunov decomposition.

Let $\alpha$  be a $\Zk$ action on a compact manifold $\M$
preserving an ergodic probability measure $\mu$. 
By commutativity, the  Lyapunov decompositions for  individual elements
of $\Zk$ can be refined to a joint invariant splitting for the action. Thus
the Multiplicative Ergodic Theorem in this case yields the following:

 \begin{proposition} 
 There are finitely many linear functionals $ \chi $ on $\Zk$, a set of
 full measure ${\cal P}$ and a splitting 
 of the tangent bundle  $T\M = \bigoplus E_{\chi}$ over ${\cal P}$  which is 
$\alpha$--invariant and measurable such that, 
  for all $a \in \Zk$ and $ v \in E_{\chi}$, the Lyapunov exponent of $v$ is $\chi 
  (a)$, ie,
 $$
   \lim _{n \rightarrow \stackrel{+}{_{-}} \infty } 
   \frac{1}{n} \log \| d a ^n  (v) \| = \chi (a),
 $$ 
  where $\| \cdot \|$ is a continuous norm on $T\M$. 
\end{proposition}

The splitting  $ \bigoplus E_{\chi} $ is called the \textit{Lyapunov decomposition\/} 
and the linear functionals $\chi$ are called the \textit{Lyapunov exponents\/} or 
\textit{Lyapunov functionals\/} of  
$\alpha$. The hyperplanes $\,\ker \,\chi \subset \Rk$ are called the \textit{Lyapunov 
hyperplanes} or \textit{Weyl chamber walls\/} and the connected components of 
$\,\Rk - \cup _{\chi} \ker \chi $ are called the \textit{Weyl chambers\/} of $\alpha$. 
The elements in the union of the Lyapunov hyperplanes are called 
\textit{singular\/}, and the elements in the union of the Weyl chambers 
are called \textit{regular\/}.

Consider  a $\Zk$ action by \textit{automorphisms\/}, of a torus or an infranilmanifold.
Then  the Lyapunov decomposition is determined by the eigenspaces of the 
automorphisms and the Lyapunov exponents are the logarithms of the moduli 
of the eigenvalues. Hence they are independent of the invariant measure, 
the Lyapunov decomposition is smooth, and the Lyapunov functionals give 
uniform estimates of expansion and contraction rates.


\subsection{Anosov $\Zk$ actions and coarse Lyapunov distributions}
 \label{Anosov actions}
  
Let $f$ be a diffeomorphism of a compact manifold $\M$. 
Recall that  $f$ is  \textit{Anosov\/} if there exist a continuous $f$--invariant
 decomposition  of the tangent bundle $T\M=E^s \oplus E^u$ and constants 
$K>0$, $\lambda>0$ such that, for all $n\in \mathbb N$,
 $$
   \begin{aligned}
  \| df^n(v) \| \,\leq\, K e^{- \lambda n} \| v \|&
      \;\text{ for all }v \in E^s, \\
  \| df^{-n}(v) \| \,\leq\, K e^{- \lambda n}\| v \|&
      \;\text{ for all }v \in E^u. 
  \end{aligned}
 $$
The distributions $E^s$ and $E^u$ are called the \textit{stable\/} and \textit{unstable\/} 
distributions of $f$. It is well-known that these distributions are tangential 
to the stable and unstable foliations $W^s$ and $W^u$ respectively.
The leaves of these foliations are $C^\infty$ injectively immersed Euclidean 
spaces. Locally, the immersions vary continuously in $\Ci$ topology. Such 
foliations are said to have \textit{uniformly $\Ci$ leaves\/}. In general, the 
distributions $E^s$ and $E^u$ are only H\"older continuous transversally 
to the corresponding foliations.

Now we consider a $\Zk$ action $\alpha$ on a compact manifold $\M$.
The action is called \textit{Anosov\/} if there is an element which acts as an 
Anosov diffeomorphism. Note that the existence of one Anosov element 
implies that all Lyapunov functionals are nonzero. Hence if such an action 
 is  \textit{algebraic\/}, ie,  by automorphisms of a torus or an infranilmanifold, 
then all regular elements are Anosov. In general, however, it is not known whether all 
regular elements of an Anosov $\Zk$ action are Anosov. To obtain a good structural 
theory for the general actions one needs to assume the existence of several 
Anosov elements.

The stable distribution of one Anosov element is invariant
under any other element, and it is natural to consider intersections of the stable
distributions for various  Anosov elements of the action. The finest such 
intersections are called \textit{coarse Lyapunov distributions\/}. This notion 
proved to be very useful for both algebraic and nonalgebraic actions.

For an algebraic action, the coarse Lyapunov distributions are defined 
everywhere and smooth. Moreover, a coarse Lyapunov distribution $E^{\chi}$ 
can be equivalently defined as a direct sum of all Lyapunov spaces
with Lyapunov functionals positively proportional to $\chi$:
\begin{equation}
 \label{(2.1)}
    E^{\chi} = \oplus E_{\chi '}, \quad \chi ' = c \, \chi \,\text{ with }\, c>0.
\end{equation}
 
For nonalgebraic actions, however, the situation is more complicated. 
It is not clear that the an intersection of several stable distributions has
constant dimension. Also, the distribution $\oplus E_{\chi '}$ in 
general is only measurable and defined almost everywhere. The next 
proposition shows that, in the presence of sufficiently many Anosov elements, 
the coarse Lyapunov distributions are well-defined, continuous, 
and tangent to foliations with smooth leaves. This is the discrete version of 
Proposition 2.4 in \cite{KaSp}. We denote the set of all Anosov elements in $\Zk$
by $\A$.

\begin{proposition} \label{CoarseLyapunov}
Let $\alpha$ be an Anosov $\Zk$--action and let $\mu$ be an ergodic  probability
measure for $\a$ with full support. Suppose that there exits an Anosov element 
in every Weyl chamber defined by $\mu$. Then for each Lyapunov exponent 
$\chi $ the coarse Lyapunov distribution can be defined as
 $$
   E^{\chi}(p) = \bigcap_{\{a \in \A | \chi (a) <0\}} \Es_a(p) = 
   \bigoplus_{\{ \chi ' = c \chi|  c>0 \}}  E_{\chi '} (p) 
  $$ 
on the set ${\cal P}$ of full measure where the Lyapunov exponents exist.
Moreover, $E^{\chi}$ is H\"{o}lder continuous, and thus it can 
be extended to a H\"{o}lder distribution tangent to the foliation
$W^\chi = \bigcap  _{\{a  \in \A \,  \mid \, \chi (a) <0\}} \Ws _a$ 
with uniformly $\Ci$ leaves.
\end{proposition}

It is easy to see that the coarse Lyapunov distributions constructed in the
proposition are indeed the finest nontrivial intersections of various stable
distributions. Extending Proposition 2.9 in \cite{KaSp}, one can also show 
that for any other invariant measure with Anosov elements in every Weyl 
chamber, \eqref{(2.1)} gives the same distributions.
Note that ergodic measures with full support always exist if the action
contains a transitive Anosov element. A natural example is given by the
measure of maximal entropy for such an element which, by uniqueness,
is invariant under the action.

Now we describe an important class of $\Zk$ actions called \textit{totally 
nonsymplectic}, or TNS. For such an action,  any pair of coarse 
Lyapunov distributions is contracted by some element. 
For an algebraic action this is equivalent to having \textit{no negatively 
proportional Lyapunov functionals}.  This property proved important and
motivated the definition for nonalgebraic actions in Katok--Ni{\c{t}}ic{\u{a}}--T{\"o}r{\"o}k \cite{KNT} similar to 
the following:

\begin{definition}
Let  $\a$ be a $\Ci$ action of $\Zk$ with a transitive  Anosov element
on a compact smooth connected manifold $\M$.  We say that $\a$ is a
 \textit{TNS Anosov action\/}  if $T\M$ splits into a direct sum of continuous 
distributions $E^i$, called the coarse Lyapunov distributions, so that:
\begin{enumerate}
    \item Any distribution $E^i$ is of the form 
    $E^i=\bigcap_{a \in \A} E_a^{\sigma (a)}$, where  $\A$ is the set of Anosov 
            elements and
            $\sigma (a)$ is either s or u.    
     \item The sum of any two distributions $E^i$ and $E^j$ is contained in 
             $E^s_a$ for some $a \in \A$.            
\end{enumerate}
\end{definition}

\begin{remark} \label{TNSremark}
Note that condition (1) ensures that the distributions 
$E^i$ are exactly the coarse Lyapunov distributions as described earlier.
It also implies that the distributions $E^i$ are H\"older continuous 
and tangent to the integral foliations $W^i$ with uniformly $\Ci$ leaves, 
called the \textit{coarse Lyapunov foliations\/}.  

It also follows that for any $E^i \not= E^j$ there exists $a \in \A$ for which
$E^i \subset E^s_a$ and $E^j \subset E^u_a$.
\end{remark}

This definition of an Anosov TNS action is very similar to the one introduced in
\cite{KNT}. In that definition, the set of all Anosov elements was replaced by a 
subcollection. It gives similar properties for the action, but the distributions 
$E^i$ may turn out to be larger than the actual coarse Lyapunov distributions. 
We consider such distributions in \fullref{Main pinching}.

\fullref{CoarseLyapunov} gives the following sufficient conditions for 
an action to be TNS Anosov. These conditions are close to being necessary.

\begin{corollary}
Let $\a$ be a $\Ci$  action of $\Zk$ with a transitive  Anosov element on a compact 
smooth connected manifold $\M$. Suppose that for some $\a$--invariant ergodic 
measure $\mu$ with full support there are no negatively proportional Lyapunov 
functionals, and in each Weyl chamber there exists an Anosov element. 
Then $\a$ is a TNS Anosov action.
\end{corollary}


\subsection{Conjugacy to algebraic models} \label{Conjugacy}

Let $f$ be an Anosov diffeomorphism of a torus or, more generally, 
an infranilmanifold $\N$. By the results of Franks \cite{Fr} and Manning \cite{M}, 
$f$ is topologically conjugate to an Anosov automorphism $A \co \N \to \N$, 
ie, there exists a homeomorphism  $\phi \co \N \to \N$ such that $A \circ \phi = \phi \circ f$. 
The conjugacy $\phi$ is unique in the homotopy class of identity. 

Now consider an Anosov $\Zk$ action $\a$ on an infranilmanifold. 
For any Anosov element of $\a$  there is a homeomorphism $\phi$ which 
conjugates it to an automorphism. It is well known that $\phi$ then conjugates 
$\a$ to an action by affine automorphisms, 
Hurder \cite[proof of Proposition 2.18]{H1}. 
This follows from
the fact that any homeomorphism commuting with an Anosov automorphism
is an affine automorphism 
(Palis and Yoccoz \cite[proof of Proposition 0]{PY}). By an affine
automorphism we mean a composition of an automorphism and a translation.
We note that an Anosov $\Zk$ action on an infranilmanifold may have no fixed
 points \cite{H1}, however, there is always a finite index subgroup which fixes
 a point and whose action is conjugate to an action by automorphisms.


\subsection{Joint integrability} \label{integrability}       

While any coarse Lyapunov distribution is integrable, a sum of two 
coarse Lyapunov distributions is not integrable in general. This can be
observed for $\Zk$ actions on nilmanifolds. However, for an algebraic 
$\Zk$ action on a torus, all Lyapunov foliations are linear, and thus 
the sum of any two coarse Lyapunov distributions is integrable. 

Let $\alpha$ be an  Anosov $\Zk$ action on a torus. Then it is topologically 
conjugate to an algebraic action (\fullref{Conjugacy}). The conjugacy maps 
the (un)stable foliations to the (un)stable foliations, and hence it maps the coarse 
Lyapunov foliations to the coarse Lyapunov foliations. We conclude that
any two coarse Lyapunov foliations for $\alpha$
are topologically jointly integrable in the following  sense. 

\begin{definition}\label{integrable} Two foliations $W^1$ and $W^2$ of a manifold 
are \textit{topologically jointly 
integrable} if there is a topological foliation $W$ such that for any $x$ the map 
 $$
    \phi \co W^1(x) \times  W^2(x) \to W(x), \quad \phi (y,z)= W^1(z) \cap W^2(y)
 $$
is a well-defined local homeomorphism. In other words, the foliations $W^1$ 
and $W^2$ give a local product structure on the leaves of $W$. 
\end{definition}

\subsection{Conformality  and uniform quasiconformality} \label{QC}  

Let $f$ be a diffeomorphism of a compact Riemannian manifold 
$\M$, and let $E$ be a continuous $f$--invariant distribution. 
The diffeomorphism $f$ is \textit{uniformly quasiconformal\/} on $E$ 
if  the quasiconformal distortion
 $$
  K_E(x,n)=\frac{\max\,\{\,\|\,df^n(v)\,\|\, :\; v\in E(x), \;\|v\|=1\,\}}
            {\,\min\,\{\,\|\,df^n(v)\,\|\, :\; v\in E(x), \;\|v\|=1\,\}}
 $$
is uniformly bounded for all $n\in\Z$ and $x\in \M$. 
We note that the notion of uniform quasiconformality does not depend on 
the choice of a Riemannian metric on the manifold.
Clearly, an Anosov diffeomorphism can be uniformly quasiconformal on $E$
only if $E$ is contained in its stable or its unstable distribution.

If $K_E(x,n)=1$ for all $x$ and $n$, the diffeomorphism is said to be 
\textit{conformal\/} on $E$.  The next result shows that any Anosov
diffeomorphism uniformly quasiconformal on $E$ is, in fact, conformal
with respect to some metric on $E$.

\begin{theorem}  \label{smooth metric} 
Let $f$ be a topologically transitive $C^\infty$ Anosov diffeomorphism  
 on a compact manifold $\,\M$. Let $E$ be a H\"older continuous 
 $f$--invariant distribution. Suppose that  $f$ is uniformly quasiconformal 
 on $E$.
 Then $f$ is \textit{conformal\/} with respect to a Riemannian metric on $E$
 which is H\"older continuous on $\,\M$.
 
 If, in addition, $E$  is tangential to a foliation $W$ with $\Ci$ leaves, 
 then this metric is  $C^\infty$ along the leaves of $W$.
\end{theorem}

The proof of this theorem is virtually identical to the proof of 
Theorem 1.3 in Sadov\-ska\-ya \cite{S}. There the stable distribution 
and the stable foliation are considered instead of $E$ and $W$.
 
Let $\Gamma$ be a group acting on $\M$ via diffeomorphisms,
and let $E$ be a continuous $\Gamma$--invariant distribution.
We say that the action is uniformly quasiconformal on $E$
if the quasiconformal distortion is uniformly bounded for all 
elements of $\Gamma$.

\begin{proposition} Suppose that the $\Gamma$--action  is generated by 
finitely many commuting diffeomorphisms. If each generator is uniformly 
quasiconformal on $E$, then the  $\Gamma$--action is uniformly 
quasiconformal on $E$.
\end{proposition}

\begin{proof}
 Any element of the action can be written as $f_1^{n_1}\dots f_k^{n_k}$, 
 where $f_1,\dots ,f_k$ are the generators. Then uniform quasiconformality
of the action follows directly from the definition.
\end{proof}


\section{Main results}  \label{Main Results}

In our first theorem we consider uniformly quasiconformal 
TNS Anosov actions.

\begin{theorem} \setobjecttype{Theorem}\label{Main QC}

Let $\a$ be a TNS Anosov action of $\Zk$ on a compact connected 
smooth manifold $\M$.
 Suppose  that: 
\begin{enumerate}
\item Any two coarse Lyapunov foliations are topologically jointly integrable.
\item The action is uniformly quasiconformal on each coarse Lyapunov distribution.
\end{enumerate}
Then a finite cover of $\a$ is $\Ci$--conjugate to 
a $\Zk$ action by affine automorphisms of a torus.

\end{theorem}

If the manifold $\M$ is a torus, then  condition (1) is automatically satisfied (see
\fullref{integrability}). Thus we obtain the following corollary. 

\begin{corollary} \label{Toral QC}
Let $\a$ be a TNS Anosov action of $\Zk$ on a torus. Suppose  that 
the action is uniformly quasiconformal on each coarse Lyapunov distribution.
Then $\a$ is $\Ci$--conjugate
 to a $\Zk$ action by affine 
automorphisms of the torus.

\end{corollary}

Note that some or all  coarse Lyapunov distributions may
be one dimensional. In this case, the quasiconformality assumption is 
trivially satisfied. In higher dimensions, quasiconformality can be replaced by
certain pinching, ie, a relationship between the slow and fast rates of
expansion/contraction. The following theorem gives a result of this type.
Note, however, that uniform quasiconformality does not 
relate the expansion/contraction rates at different points,
so it does not imply the 1/2--pinching below.

\begin{definition} Let $f$ a diffeomorphism which contracts an
invariant distribution $E$.
We say that $f$ is \textit{1/2--pinched\/} on  $E$ if there exist numbers $K>0$ and
$0<\mu< \lambda\,$ with $\,\lambda <2\mu\,$ such that for any $v\in E$,
$$ 
   K^{-1} e^{-n\lambda} \|v\| \,\le\,    \| df^n(v) \| \,\le\, K e^{-n\mu} \|v\|.
$$  
\end{definition}

\begin{theorem}\setobjecttype{Theorem} \label{Main pinching}
Let $\a$ be an Anosov action of $\Zk$ on a compact connected smooth 
manifold $\M$. Suppose  that  $T\M$ splits into a direct sum of continuous 
distributions $E^i$, where each $E^i$ is an intersection of stable  distributions 
of some Anosov elements of the action, 
and that for any two distributions $E^i$ and $E^j$:
\begin{enumerate}
\item The corresponding foliations  are topologically jointly integrable.
\item There is  an element which contracts both distributions and is 
         1/2--pinched on $E^i$.
\item  There is an element $\,a \in \Zk$  which expands $E^i$ faster than 
          it expands $E^j$, \\ ie, for any $x \in \M$,
   $$   
          \max  \,\{ \| da(v)\| :\; v\in E^j(x), \;\|v\|=1 \} \, < \,
           \min   \,\{ \|da(v)\| :\; v\in E^i(x), \;\|v\|=1 \}.
   $$
\end{enumerate}
Then $\a$ is $\Ci$--conjugate to a $\Zk$ action by affine automorphisms 
of an infranilmanifold.
\end{theorem}

Applying this theorem to actions on tori, we obtain the following result.
 
\begin{corollary} \label{Toral pinching}
   Let $\a$ be an Anosov action of $\Zk$ on a torus satisfying conditions 
   (2) and (3) of \fullref{Main pinching}. Then $\a$ is $\Ci$--conjugate to 
   a $\Zk$ action by affine automorphisms of the torus.
\end{corollary}

\begin{remark} Condition (3) is used in \fullref{Main pinching} only 
to obtain certain smoothness of $E^i$ along the leaves of the foliation $W^j$ 
tangential to $E^j$. It can be substituted by the assumption that the regularity 
of $E^i$ along $W^j$ is $C^{1,\beta}$ for all $\beta<1$.
\end{remark}

In \fullref{Main pinching}, each distribution $E^i$ may be a direct sum 
of several coarse Lyapunov distributions (see the discussion after \fullref{TNSremark}). If the distributions $E^i$ are the coarse Lyapunov 
distributions, then the infranilmanifold is finitely covered by a torus.

\begin{corollary}\setobjecttype{Corollary} \label{TNS pinching}
Let $\a$ be a TNS Anosov action of $\Zk$ on a compact connected 
smooth manifold $\M$. Suppose  that the coarse Lyapunov distributions satisfy
conditions (1), (2), and (3) of \fullref{Main pinching}. 
Then a finite cover of $\a$ is $\Ci$--conjugate to 
a $\Zk$ action by affine automorphisms of a torus.

\end{corollary}

We can also deduce some results on local rigidity of algebraic Anosov actions.
The local rigidity of higher rank Anosov $\Zk$ actions by automorphisms 
of tori was proved by Katok and Spatzier \cite{KS97} under the assumption that 
the automorphisms are semisimple (ie, have no Jordan blocks). The general case 
has not been resolved. In a recent preprint \cite{EF}, Einsiedler and Fisher
considered the Jordan block case under the 1/2--pinching assumption on each
coarse Lyapunov distribution.  \fullref{perturbation}  gives an 
alternative proof of their result in the TNS case. 

\begin{corollary}\setobjecttype{Corollary} \label{perturbation}
Let $\a$ be a TNS Anosov action of $\Zk$ by toral automorphisms.
Suppose that for each  coarse Lyapunov distribution there exists a regular
element which is 1/2--pinched on this distribution. 

Then $\a$ is locally rigid, ie, it is $\Ci$--conjugate to any  sufficiently 
$C^1$--small perturbation.
\end{corollary}

More generally,   \fullref{TNS pinching} can be applied 
to actions which have pinching similar to that of small perturbations
of 1/2--pinched algebraic actions.

\begin{corollary}\setobjecttype{Corollary} \label{functional pinching}
Let $\a$ be a TNS Anosov action of $\Zk$ on a compact connected 
smooth manifold $\M$. Suppose  that: 
\begin{enumerate}
  \item Any two coarse Lyapunov foliations are topologically jointly integrable.
  \item For any coarse Lyapunov distribution $E$ there exist  functionals
      $\chi_{\min}$ and $\chi_{\max}$ proportional with a constant $1\le c <2$,
      such that for some $K>0$
      $$
        K^{-1} e^{\chi_{\min} (a)} \|v\| \,\le\, 
        \| da (v) \| \,\le\, K e^{\chi_{\max} (a)} \|v\|
     $$ 
for any $a \in \Zk$ with $\chi_{\min}(a) >0$ and any $v \in E$.
\end{enumerate}
Then a finite cover of $\a$ is $\Ci$--conjugate to 
a $\Zk$ action by affine automorphisms of a torus.

\end{corollary}

\begin{remark}
The inequality in condition (2) of \fullref{functional pinching} of  can be replaced by a weaker assumption
$$
 K^{-1} e^{\chi_{\min} (a) - \e |a|} \|v\| \,\le\, 
 \| da (v) \| \,\le\, K e^{\chi_{\max} (a) + \e |a|} \|v\| 
$$ 
if $\e$ is small enough for the given system of functionals.
\end{remark}


\section[Proof of \ref{Main QC}]{Proof of \fullref{Main QC}}\label{Proof Main QC}

\subsection{Outline  of the proof}  

We begin by considering any two coarse Lyapunov distributions 
$E^1$ and $E^2$, and the corresponding coarse Lyapunov foliations 
$W^1$ and $W^2$. In Sections \ref{regularity of sum} and 
\ref{sec smooth distributions} we show that the leaves of 
$W=W^1 \oplus W^2$ are smooth, and that $E^1$ and $E^2$
are $C^{1,\beta}$ along these leaves.

In \fullref{linearize} we discuss the nonstationary linearization of the
action along the leaves of $W^1$ and the associated affine structures
on these leaves. In the next section we show that the holonomies
of the foliation $W^2$ between the leaves of $W^1$ are affine.
This implies that $E^2$ is $\Ci$ along $W^1$. It follows that 
all coarse Lyapunov distributions are $\Ci$ on $\M$.

In \fullref{smooth linearizations} we show that the nonstationary 
linearizations depend smoothly on the base point. We use this to construct
a $\Ci$ invariant affine connection in \fullref{connection}. Together
with the smoothness of the coarse Lyapunov distributions, this implies that
the action is $\Ci$--conjugate to an action by affine automorphisms of 
an infranilmanifold. Finally, we show that this infranilmanifold is finitely
covered by a torus.


\subsection{Regularity of the sum of two coarse Lyapunov foliations} 
\label{regularity of sum}

Since $\a$ is a TNS Anosov action, the coarse Lyapunov distributions 
are defined everywhere on the manifold and H\"older continuous.
They are tangential to coarse Lyapunov foliations with uniformly $\Ci$ 
leaves. (See \fullref{Anosov actions}.)

We consider any two coarse Lyapunov distributions $E^1$ and $E^2$, 
and the corresponding coarse Lyapunov foliations $W^1$ and $W^2$. 
By assumption (1) of the theorem, the foliations  $W^1$ and $W^2$ 
are topologically jointly integrable, see \fullref{integrable}.
We will  write $W^1 \oplus W^2$ for the foliation $W$ in the definition. 
Since the action is TNS, there is an element 
which contracts  $W$. Applying the iterates of its inverse, one can see
that the local product structure on the leaves of $W$ extends to the 
global one. 

\begin{lemma} The topological foliation $W=W^1 \oplus W^2$
has uniformly $\Ci$ leaves tangent to the distribution $E=E^1\oplus E^2$.
\end{lemma}

\begin{proof}
This follows from  the Journ\'e lemma below.
A function is said to be \textit{uniformly $\Ci$\/} along a foliation with smooth
leaves if the partial derivatives of all orders in the foliation directions
exist and are continuous on the manifold. 

\begin{lemma}{\rm \cite {J}}\qua  \label{Journe} Let $\M$ be a $\Ci$ manifold, 
and let $W^1$ and $W^2$ be  two transverse H\"older continuous 
foliations with uniformly $\Ci$ leaves. If a function $\phi$ is uniformly 
$\Ci$ along the leaves of the two foliations, then it is $\Ci$ on $\M$.
\end{lemma}

We recall that the leaves of $W^1$ and $W^2$ are $\Ci$ immersed manifolds 
tangent to the distributions $E^1$ and $E^2$ respectively. Moreover, they are 
uniformly $\Ci$. 
We can identify a small neighborhood of a point $x \in \M$
with a ball in $\Rn$ and consider the orthogonal projection $\Pi$ to the 
subspace $E(x)=E^1(x) \oplus E^2(x)$. The restriction of $\Pi$ to the 
leaf $W(x)$ is locally a homeomorphism that projects foliations $W^1$ 
and $W^2$ to transverse foliations $\wwbar W^1$ and $\wwbar W^2$ in 
$E(x)$ with uniformly $\Ci$ leaves.  Then the inverse $\phi$ of this
restriction is uniformly $\Ci$ along the leaves of $\wwbar W^1$ and 
$\wwbar W^2$ and hence it is $\Ci$ by the Journe lemma.
This implies that $W(x)$ is a $\Ci$ immersed manifold. In particular,
its tangent distribution is $E=E^1 \oplus E^2$. Since the immersions of 
$W^1(x)$ and $W^2(x)$ depend continuously on $x$ in $\Ci$ topology
it is easy to see that so do the immersions of $W(x)$.
\end{proof}


\subsection{Coarse Lyapunov distributions are $C^{1,\beta}$ along $W$}
\label{sec smooth distributions} 

Let $E^1$ and $E^2$ be two coarse Lyapunov distributions,
let $E=E^1 \oplus E^2$ be their sum, and  let $W=W^1 \oplus W^2$  be the  
corresponding  foliation. Above we established that the leaves of this foliation 
are $\Ci$ manifolds. We study regularity of the distributions $E^1$ and 
$E^2$ along the leaves of $W$.

\begin{proposition} \label{smooth distributions} 
Suppose that there exists an element $f$ that expands $E^1$, 
contracts $E^2$ and is uniformly quasiconformal on $W^2$.
Then $E^1$ is $C^{1,\beta}$  along the 
leaves of $W$ for some $\beta>0$. More precisely,
 its first derivatives in $W$ 
directions exist and are H\"older continuous along the leaves 
of $W$ with  uniform bounds on $\M$ for the derivatives and 
the H\"older constants.
\end{proposition}

\begin{proof}
We use the $C^r$ Section Theorem of M Hirsch, C Pugh, and M Shub. 
See Theorem 3.2, as well as its more detailed version Theorem 3.5,  
and Remarks 1 and 2 after Theorem 3.8 in  \cite{HPS}.  

We consider the distribution  $E=E^1\oplus E^2$. This is a continuous bundle
over $\M$ which is $\Ci$ smooth along the leaves of $W$. There exist distributions
$\bar E^1$ and $\bar E^2$ in $E$ which are close to $E^1$ and $E^2$
respectively and  $\Ci$ along the leaves of $W$. Now we can consider a vector
bundle $\L$ whose fiber over $x$ is the set of linear operators from $\bar E^1(x)$ 
to  $\bar E^2(x)$. The differential of $f$ induces a natural action $F$ on $\L$.   
Suppose that $F$ contracts fibers of $\L$, ie, for any $x \in \M$ and any $u,v \in \L(x)$ 
 $$
   \|F(u)-F(v)\|_{fx} \le k_x \|u-v\|_x \;\text{ with } \sup_{x \in M}  k_x <1
 $$
with respect to some norm on $\L$.
Then Theorem 3.1 in \cite{HPS} gives  existence, uniqueness, 
and continuity of  the invariant bounded section. 
By uniqueness, the graph of this invariant section is the distribution $E^1$.
Moreover, the $C^r$ Section Theorem \cite[Theorems 3.2 and 3.5]{HPS}
states that if $\L$ is $C^r$ and 
\begin{equation}
   \sup_{x \in M}  k_x \a_x ^r <1 \;\text{ where }
    \a_x = \|(df_x)^{-1}\|, \label{(4.1)}
\end{equation}
then the invariant section is $C^r$ smooth. Theorem 3.8 in \cite{HPS} and Remarks
1 and 2 after it imply that this conclusion holds for noninteger values of $r$ in the
H\"older category, and that the compactness of the base can be replaced by
the boundedness condition as in Theorem 3.2 in \cite{HPS}. 

We are interested only in the smoothness of the invariant section along 
the leaves of $W$.
In this case, we need the smoothness of the bundle $\L$ only along $W$ 
and we can use $\a_x = \|(df|_{E(x)})^{-1}\|$ in \eqref{(4.1)}. 
This can be also seen by formally
applying the $C^r$ Section Theorem to the bundle $\L$ over  the
disjoint union of the leaves of $W$ as the base.

Now we use the assumption that $f$ is uniformly quasiconformal 
on $W^2$ to verify \eqref{(4.1)}.
We fix some continuous norm on the distribution $E$ and endow the 
fibers of $\L$ with the standard operator norm. We denote
\begin{gather*}
   \,\, l_x = \min \,\{\,\|df_x (v) \| : \; v \in E^1(x),\;  \|v\|=1\,\}\\
  m_x = \min \,\{\,\|df_x (v) \| : \; v \in E^2(x), \; \|v\|=1\,\}\\
  M_x = \max \,\{\,\|df_x (v) \| : \; v \in E^2(x),\;  \|v\|=1\,\}.
\end{gather*}
$$
     \a_x = \|(df|_{E(x)})^{-1}\|=1/m_x \quad \text{and} \quad k_x \approx M_x/l_x
\leqno{\hbox{Then}} 
    $$
since $\bar E^1$ and $\bar E^2$ are close to $E^1$ and $E^2$.
Hence
$$
    k_x \a_x ^{1+\beta} \, \approx \,\frac{M_x}{l_x m_x^{1+\beta}} \, \le
    \, \frac{1}{l\cdot m_x^{\beta}} 
    \cdot \frac{M_x}{m_x} \,\le\, \frac{K}{l \cdot (\inf_x \{m_x\})^{\beta}}
 $$
where $l=\inf_x l_x$ and 
 $$
 K= \sup_{x, \,n}\,\frac{\max\,\{\,\|\,df^n(v)\,\|\, :\; v\in E^2(x), \;\|v\|=1\,\}}
            {\,\min\,\{\,\|\,df^n(v)\,\|\, :\; v\in E^2(x), \;\|v\|=1\,\}}
$$ 
is the quasiconformal distortion bound. 
Note that K is a uniform bound for all iterates $f^n$. Since $f$ expands 
$E^1$ we can replace it by  $f^n$, if necessary, to ensure that $l$ is large
enough so that $K/l <1$. Then the right hand side
is less than 1 for some $\beta >0$. Once this iterate is chosen, we can
take $\bar E^1$ and $\bar E^2$ close enough to $E^1$ and $E^2$ to
guarantee that $ \sup_{x \in M}  k_x \a_x ^{1+\beta} <1$. Hence, by the 
$C^r$ Section Theorem, the distribution $E^1$ is $C^{1,\beta}$ smooth
along the leaves of $W$.
\end{proof}

By \fullref{TNSremark}, for any two coarse Lyapunov distributions
 there exists an element of the action 
which expands the first distribution and contracts the second one.
Thus,  under the assumptions of 
\fullref{Main QC}, \fullref{smooth distributions} implies that 
any two coarse Lyapunov
foliations  $W^1$ and $W^2$ are $C^{1,\beta}$ smooth along the leaves 
of their sum.


\subsection{Nonstationary linearizations} \label{linearize}

For each coarse Lyapunov foliation $W=W^i$ we use the following proposition 
to obtain a nonstationary linearization $h=h^i$ of the action
along the leaves of $W$. 

\begin{proposition}{\rm \cite[Proposition 4.1]{S}}\qua \label{linearization}
Let $f$ be a diffeomorphism of a compact Riemannian manifold $\M$. 
Let $W$ be a continuous invariant foliation with $C^\infty$
leaves, and let $E$ be its tangent distribution. 

Suppose that $\|\,df|_E\,\|<1$, and there exist
$K>0$ and $\e >0$ such that for any $x\in \M$ and $n\in\mathbb N$, 
\begin{equation}
  \|\left( df^n|_{E(x)} \right)^{-1}\| \cdot \|\,df^n|_{E(x)}\|^2 
  \leq K(1-\e )^n. \label{(4.2)}
\end{equation}
Then for any $x\in \M$ there exists a $\Ci$ diffeomorphism 
$h_x\co W(x) \to E(x)$ such that
  \begin{enumerate}
    \item  $h_{fx}\circ f=df_x \circ h_x$;
    \item $h_x(x)=0$ and $dh_x (x)$ is the identity map; 
    \item $h_x$ depends continuously on $x$ in $\Ci$ topology.
 \end{enumerate}
\end{proposition}
 
 Let $f$ be an element of the action  which contracts a coarse 
 Lyapunov foliation $W$. 
Since $f$ is uniformly quasiconformal on $W$, 
 $$
   \|\left( df^n|_{E(x)} \right)^{-1}\| \cdot \|\,df^n|_{E(x)}\| 
 $$ 
is uniformly  bounded in $x$ and $n$. Hence \eqref{(4.2)} is satisfied and there exists a linearization $h$ 
for $f$ along the leaves of $W$.
Since such a linearization is unique \cite[Lemma 4.1]{S}, 
$h$ linearizes any diffeomorphism  which commutes with $f$. 
Indeed, if $g\circ f=f\circ g$, then it is easy to see that 
$dg^{-1} \circ h \circ g\,$ also gives a 
linearization for $f$, and hence $\,h \circ g=dg \circ h$.
Therefore,  $h$ provides linearization for all elements of the action.

\begin{lemma} \label{ratio}
Under the assumptions of \fullref{linearization},
for any $R>0$ there exists $K>0$ such that
for any two points $x$ and $z$ on the same leaf of $W$
 with $\dist\,(x,z) <R$ and any $n > 0$
$$
   \| df^n|_{E(z)} \| \le K \cdot \|df^n|_{E(x)} \|.
$$
\end{lemma}
\begin{proof}  Using the linearization $h$ along the leaves of $W$ we can write 
\begin{equation*}
\begin{split}
    f^n|_{W(x)}&=(h_{f^n x})^{-1} \circ df^n|_{E(x)} \circ h_x,\\
 \text{so} \qua df^n|_{E(z)}&=( dh_{f^n x}(f^nz) )^{-1} \circ df^n|_{E(x)} \circ dh_x (z).
\end{split}
\end{equation*}
Since $h_x$ depends continuously on $x$ in $\Ci$ topology, $\|dh_x (y)\|$
and $\|(dh_x (y))^{-1}\|$ are uniformly bounded for all $x \in \M$ and $y \in W(x)$
with $\dist\,(x,y) <R$. Since $\dist\,(f^n x,f^n z) <\dist\,(x,z) <R$, the norms
of the first and last term in the right hand side are uniformly bounded and the
lemma follows.
   \end{proof}

\begin{proposition} \label{consistency}
Under the assumptions of \fullref{linearization}, the map
$$
    h_z \circ  h_x^{-1}\co E(x) \to E(z)
$$
is affine for any $x$ and $z$ on the same leaf of $W$. Hence the 
nonstationary linearization $h$ defines affine 
structures on the leaves of $W$. 
\end{proposition}

By an affine structure we understand an atlas with affine transition maps.

\begin{proof}
It suffices to show that 
 $$
   d\left( h_z \circ  h_x^{-1}\right) (\bar y)=d \left( h_z \circ  h_x^{-1} \right) (0_x)\quad
  \text{for any }\bar y \in E(x).
 $$
 Since $h_{x}=\left( df^n|_{E(x)} \right) ^{-1} \circ h_{f^n x} \circ f^n|_{W(x)}$ we obtain 
 $$
    h_z \circ  h_x^{-1} =\left( df^n |_{E(z)} \right)^{-1} \circ h_{f^n z} \circ 
     (h_{f_nx})^{-1} \circ df^n |_{E(x)}.
 $$  
Then for any $\bar y \in E(x)$
 $$
    d\left( h_z \circ  h_x^{-1} \right) (\bar y) = \left( df^n |_{E(z)} \right)^{-1} \circ 
    dh_{f^n z} (f^n y) \circ \left( dh_{f_nx}(f^n y) \right) ^{-1} \circ df^n |_{E(x)},
 $$  
where $y =(h_x)^{-1}(\bar y)\in  W(x)$. Hence
 $$
 \begin{aligned}
     \|\, d ( h_z \circ  h_x^{-1} ) (\bar y) - d(  h  _z \circ h_x^{-1} ) (0_x) \,\| &\le \\
    \|\,\left( df^n |_{E(z)} \right)^{-1}\,\| \cdot 
        \|\, dh & _{f^n z} (f^n y) \circ (dh_{f_nx}(f^n y))^{-1} - \\
        dh &_{f^n z} (f^n x) \circ ( dh_{f_nx}(f^n x) ) ^{-1} \,\| 
    \cdot \|\, df^n |_{E(x)} \|.
  \end{aligned}
 $$ 
Note that all four differentials of $h$ in the middle term are close to the identity $\Id$
(in fact $dh_{f_nx}(f^n x)=\Id$). More precisely, $dh_x(y)$ depends 
Lipschitz continuously on $y \in W(x)$  with $\dist\,(x,y) <R$ and the Lipschitz 
constant is uniform in $x \in \M$. Hence the norm of the difference between
each of these four differentials and $\Id$ is of order 
 $$
    \max \{ \dist\,(f^n x,f^n y), \dist\,(f^n y,f^n z) \} \le K_1 \|\, df^n |_{E(z)}\|,
 $$
as easily follows from \fullref{ratio}. Applying this lemma again we obtain
\begin{multline*}
     \| d(h_z \circ  h_x^{-1})(\bar y) - d ( h_z \circ   h_x^{-1})(0_x) \| \le \\
     \| (df^n |_{E(z)})^{-1}\| \cdot  K _2 \|  df^n  |_{E(z)}  \|  \cdot K  \| df^n |_{E(z)}\| 
     \to 0
\end{multline*} 
 by \eqref{(4.2)}. This shows that the differential of $h_z \circ  h_x^{-1}$ 
 is constant on $E(x)$ and thus the map is affine.
\end{proof}

We conclude that for every coarse Lyapunov foliation there exist affine
structures on its leaves. The maps induced by the elements of the action 
on the leaves are affine with respect to these structures.


\subsection{Holonomies are affine and 
                   coarse Lyapunov distributions are $\Ci$} \label{holonomies}  

We continue to study the regularity of two coarse Lyapunov foliations
$W^1$ and $W^2$ along the leaves of their sum $W=W^1 \oplus W^2$. 
In  this section we consider holonomy maps $H$ between the leaves 
of $W^1$ along the leaves of $W^2$. Let  $x$ and $z$ be  two nearby points 
on the same leaf of $W^2$. For a point $y$ in $W^1(x)$ close to $x$,
we denote
 $$
 H(y)=H_{xz}(y)=W^2(y) \cap W^1(z).
 $$   
Our goal is to show that the holonomies are affine with respect to the affine 
structures on the leaves of $W^1$. 

\begin{proposition} \label{affine holonomies}
Let $f$ be a diffeomorphism which contracts $W^1 \oplus W^2$. 
Suppose that there exists $0 < \beta \le 1$ such that the holonomy maps are 
$C^{1,\beta}$ and for any two nearby points $x$ and $z$ on the same leaf of $W^2$,
\begin{equation}
    \|\left( df^n|_{E^1(z)} \right)^{-1}\| \cdot \|\,df^n|_{E^1(x)}\|^{1+\beta} \to 0 \quad 
    \text{as }\, n\to \infty. \label{(4.3)} 
\end{equation}
Then the holonomy maps are affine with respect to the affine structures 
on the leaves of $W^1$ and are uniformly $\Ci$, ie, they depend continuously 
on $x$ and $z$ in $\Ci$ topology. Also, the distribution $E^2$ is uniformly $\Ci$
along the leaves of $W$.
\end{proposition}

\begin{proof} Since the diffeomorphism $f$  contracts both
 $W^1$ and $W^2$, the holonomy map 
$H=H_{xz} \co W^1 (x) \to W^1 (z)$ 
is defined on the whole leaf $W^1(x)$.
Indeed, for any $y$ in $W^1(x)$ we have
$H(y)=\left( f^{-n} \circ H_{f^n x \, f^n z} \circ f^n \right) (y)$, where $n$ is 
large enough so that the points $f^nx$, $f^n y$, and $f^n z$ are close.

Let $h_x\co W^1(x) \to E^1(x)$ and $h_{z}\co W^1(z) \to E^1(z)$ be 
the linearizations given by \fullref{linearization}. 
Consider  the map 
 $$
    \bH=\bH_{xz}=h_z \circ H \circ (h_x)^{-1}\co E^1(x) \to E^1(z).
 $$
 Our goal is to show that $\bH$ is linear and hence $H$ is affine.
 Note that $\,\bH(0_x)=0_z$ since $h_x(x)=0_x$ and $h_z (z)=0_z$.
 For a point  $y$ in $W^1(x)$, we denote $\bar y=h_x(y)$. 
To prove that the map $\bH$ is linear it suffices to show that 
$\,d \bH (\bar y) =d \bH (0_x)$ for any $y$ in $W^1(x)$. 

Using forward iterations of the  diffeomorphism $f$, we can write
\begin{gather*}
     \bH= \left( df^n|_{E^1(z)} \right) ^{-1} \circ \bH_n \circ \, df^n|_{E^1(x)},\\
     \bH_n= h_{f^n z} \circ H_n \circ (h_{f^n x})^{-1}\co E^1(f^n x) \to E^1(f^n z),\tag*{\hbox{where}}
\end{gather*}
and $H_n\co W^1(f^n x) \to W^1 (f^n z)$ is the holonomy map along
the leaves of $W^2$. Then 
  $$
     d\bH= \left( df^n|_{E^1(z)} \right) ^{-1} \circ d\bH_n \circ \, df^n|_{E^1(x)}.
  $$
In $E^1(f^n x)$ we denote $0_n =0_{f^n x}$ and 
$ \bar y_n=df^n_x(\bar y)=(h_{f^n_x})^{-1}(f^n y)$.   Then 
  \begin{gather*}
    d \bH (\bar y)=\left( df^n|_{E^1(z)} \right)^{-1} \circ d\bH_n (\bar y_n) 
    \circ \, df^n|_{E^1(x)}\\
   d \bH(0_x)=\left( df^n|_{E^1(z)} \right) ^{-1} \circ d\bH_n (0_n)  
    \circ \, df^n|_{E^1(x)}. \tag*{\hbox{and}}
  \end{gather*}
We estimate the norm of the  difference between 
$d \bH (\bar y)$ and $d \bH (0_x)$ as follows.
 $$
\begin{aligned}
  \| d \bH (\bar y)-d \bH (0_x)\| = 
  \| & \left( df^n|_{E^1(z)} \right) ^{-1} \circ  \left(d\bH_n (\bar y_n) - d\bH_n (0_n) \right) 
        \circ \, df^n|_{E^1(x)} \| \\
 \leq  \| & \left( df^n|_{E^1(z)} \right) ^{-1} \| \cdot \|\, d\bH_n (\bar y_n) - d\bH_n (0_n) \, \| 
        \cdot  \| \, df^n|_{E^1(x)} \| 
\end{aligned}   
 $$   
Since the holonomy maps are $C^{1,\beta}$ with uniform H\"older constant,
so are $\bH$. Hence
  $$
     \| \, d\bH_n (\bar y_n) - d\bH_n (0_n) \, \| \,\le\, 
          K \cdot ( \dist(\bar y_n, 0_n) )^\beta \le 
      K \cdot  (\dist(\bar y, 0_x))^\beta \cdot \|\, df^n|_{E^1(x)} \|^\beta.
 $$
Therefore 
$$
\|\, d \bH (\bar y)-d \bH (0_x) \,\| \le 
  K \cdot  (\dist(\bar y, 0_x))^\beta \cdot \| \left( df^n|_{E^1(z)} \right) ^{-1} \|  
          \cdot \| \, df^n|_{E^1(x)} \|^{1+\beta} \to 0
$$
as $n \to \infty$.
This implies that $d \bH (\bar y)=d \bH (0_x)$ for any $\bar y$ in $E^1(x)$. 
This shows that $d\bH$ is constant on  $E^1(x)$, and therefore 
the map $\bH$ is linear. 

Since $h_x$ and  $h_z$ are $\Ci$ diffeomorphisms, 
the holonomy map $H= (h_z)^{-1}\circ \bH \circ h_x \co  W^1(x) \to W^1(z)$ 
is also a $\Ci$ diffeomorphism.
Recall that the linearization depends continuously on the base point in $\Ci$ 
topology. Then, since $\bH$ is linear and depends continuously on $x$ and $z$,
the holonomy map $H=H_{xz}$ depends continuously on $x$ and $z$ in 
$\Ci$ topology. This implies that $E^1$ is uniformly $\Ci$
along the leaves of $W$ \cite{PSW}.
\end{proof}

\begin{corollary} \label{smooth E}
Under the assumptions of \fullref{Main QC}, all coarse Lyapunov 
distributions $W^i$ are $\Ci$ on $\M$. Also, the holonomy maps 
between the leaves of $W^i$ along the leaves of $W^j$ are affine.
\end{corollary}

\begin{proof}
We fix coarse Lyapunov foliations $W^i$ and $W^j$.
To apply \fullref{affine holonomies} we consider an element $f$ 
of the action which contracts $E^i \oplus E^j$. Such an element exists since the 
action is TNS, and it is uniformly quasiconformal on $E^i$ by assumption 
(2) of the theorem.  In \fullref{sec smooth distributions} we 
showed that $E^j$ is $C^{1,\beta}$ along the leaves of $W=W^i \oplus W^j$
with some $\beta  >0$ and a uniform H\"older constant on $\M$. Hence so 
are the holonomy maps $H$ along $W^j$ \cite{PSW}. 
By \fullref{df^n conf}, $\| \left( df^n|_{E^i(z)} \right)^{-1} \| \cdot  \|\, df^n|_{E^i(x)} \|$ 
is bounded by a constant independent of $n$, and hence 
  $$
     \| \left( df^n|_{E^i(z)} \right)^{-1} \| \cdot  \|\, df^n|_{E^i(x)} \|^{1+\beta} \to 0
  $$ 
for any $\beta >0$. 
Now it follows from \fullref{affine holonomies} that the holonomy maps 
are affine and uniformly $\Ci$, and  $E^j$ is uniformly $\Ci$ along the leaves of $W$.

In particular, we have established that any coarse Lyapunov
distribution is uniformly $\Ci$ along the leaves of any coarse
Lyapunov foliation.  To conclude that any coarse Lyapunov distribution
is $\Ci$ on $\M$ we use the following well-known result (see Katok and
Lewis \cite{KL2} and Goetze and Spatzier \cite{GS}), which is obtained
by inductive application of Journ\'e \fullref{Journe}.

\begin{lemma} \label{Ind Journe} 
Let $\phi$ be a map from $\M$ to a finite-dimensional $\Ci$ manifold.
If $\phi$ is uniformly $\Ci$ along the leaves of every coarse Lyapunov 
foliation then it is $\Ci$ on $\M$.
\end{lemma}

To complete the proof of \fullref{smooth E}, it remains to establish 
the following lemma.

\begin{lemma} \label{df^n conf} Let $E$ be a H\"older continuous  distribution 
invariant under a diffeomorphism $f$. Suppose that  $f$ is uniformly 
quasiconformal on  $E$. Let $x$ and $z$ be two  points such that 
for some $K_1$ and $0<\la <1$, $\;\dist \, (f^n x, f^n z) \le K_1 \cdot \la ^n$ for all $n \ge 0$.
Then  there exists a constant $K$ such that for all $n \ge 0$,
\begin{equation}
  \| \left( df^n|_{E(z)} \right)^{-1} \| \cdot  \|\, df^n|_{E(x)} \| \le K. \label{(4.4)}
\end{equation}
\end{lemma}

Note that the estimate \eqref{(4.4)} implies the one in \fullref{ratio}, and in the 
uniformly quasiconformal case the two estimates are, in fact, equivalent.
The main difference is that in \fullref{df^n conf} the points $x$ and $z$ 
are not required to be on the same leaf of a foliation tangential to $E$.
Hence the proof is different and requires a stronger assumption.

\begin{proof}
Since $f$ is uniformly quasiconformal on $E$, \fullref{smooth metric} 
implies that it is conformal with respect to a H\"older continuous Riemannian
metric on $E$. We will use this metric in the proof. Note that (4.4) is independent 
of the choice of a continuous metric on $E$.

To simplify the notations, in this proof only we write 
$df^i_x$ instead of $df^i|_{E(x)}$.
Since $E$ is H\"older continuous, $df|_{E}$ is also H\"older continuous with 
some exponent $\b >0$. Thus we obtain
 $$
  \frac{\|df_x\|}{\|df_z\|} \;\le\; 1+ \frac{|\, \|df_x\|-\|df_z\| \,|}{\|df_z\|} \;\le\;
    1+ \frac{K_2 \cdot (\dist (x,z))^\beta}{\min_z \|df_z\|} = 
  1+ K_3 \cdot (\dist (x,z))^\beta.
 $$ 
Since $f$ is conformal on $E$, $\,\|(df_z)^{-1} \| = \|df_z\|^{-1}$, 
and we can estimate
$$\eqalignbot{
\| ( df^n_z& )^{-1} \| \cdot  \|\, df^n_x\,\| \cr
  &\le \|(df_z)^{-1}\| \cdot \|(df_{fz})^{-1}\|  \dots 
     \|(df_{f^{n-1}z})^{-1}\| \cdot
      \|df_x\| \cdot \|df_{fx}\| \dots  \|df_{f^{n-1}x} \| \cr
 &= \frac{\|df_x\|}{\|df_z\|} \cdot  \frac{\|df_{fx}\|}{\|df_{fz}\|} 
       \cdot\cdot\cdot \frac{\|df_{f^{n-1}x} \|}{\|df_{f^{n-1}z} \|}  \le
      \prod_{i=0}^{n-1} 
      \left( 1+K_3 \cdot \left( \dist (f^i x, f^i z) \right)^\b \right) \cr
  &\le  \prod_{i=0}^{n-1} 
      \left( 1+K_3 \cdot \left( K_1  \cdot \la ^i \right)^\b 
      \right)  \le   
      \prod_{i=0}^{n-1} \left( 1+K_4 \cdot \la^{i\b} \right)  \le    K.\cr}
\proved
$$
\end{proof}


\subsection{The linearizations depend smoothly on the basepoint} 
\label{smooth linearizations} 

Let $W^1$ and $W^2$ be two coarse Lyapunov foliations 
and let $W=W^1 \oplus W^2$ be their sum.
In \fullref{holonomies} we established that the holonomy maps 
between the leaves of $W^1$ along $W^2$ are affine.
We will use this fact to show that the linearizations $h^1_x: W^1(x)\to E^1(x)$ 
depend smoothly on the base point along the leaves of $W$.

\begin{proposition} \label{smooth h}
If $W^1$ and $W^2$ are two coarse Lyapunov foliations
then the linearizations $h^1_x$ depend uniformly $\Ci$
on the base point along the leaves of $W^1 \oplus W^2$.
\end{proposition}

\begin{proof} First we construct a local linearization $h$ along $W$ using linearizations
$h^1$  along $W^1$ and $h^2$ along $W^2$. 
Let $E$ be the distribution tangent to $W$ and 
let $U$ be a small open  neighborhood of a point $x$ in $W$.
We define the map $h_x\co U \to E(x)$ as follows: 
 $$ 
  h_x|_{W^1(x)}=h^1_x, \quad h_x|_{W^2(x)}=h^2_x,  
 $$
and for a point $p$ in $U$ we set
 $$
  h_x(p)=h^1_x (y)+ h^2_x(z),  \text{ where }
  y=W^2(p) \cap  W^1(x), \; z= W^1(p) \cap  W^2(x). 
 $$
 Since the foliations $W^1$ and $W^2$ are $\Ci$ smooth along the leaves 
 of $W$, the map $h_x$ is a $\Ci$ diffeomorphism.
 It is easy to see that the family of maps $h_x$
satisfies conditions (1), (2), and (3) of \fullref{linearization}
and thus $h$ gives a local nonstationary linearization along $W$.
However, we do not use properties (1) and (2) in the proof.

Let us identify $E(x)$ with 
$\Rn \times \Rm$ in such a way that $E^1(x)$ corresponds to $\Rn$, and 
$E^2(x)$ corresponds to $\Rm$. Then $h_x$ identifies the neighborhood 
$U\subset W$ of $x$ with an open neighborhood $\bar U$ of 0  
in $\Rn\times\Rm$. It is clear  from the construction of $h_x$ that
the leaves of $W^1$ and $W^2$ correspond to subspaces parallel to 
$\Rn$ and $\Rm$ respectively. For a point $p$ in $U$, the tangent space $E_p$
is  identified with $\Rn\times\Rm$ by $(dh_x)_p$
in such a way that $0_p$ corresponds to $\bar p = h_x(p)$.
We will show that the maps $h^1_p$, $\,p\in U$, are identity maps when 
written  in these local coordinates. In other words, this coordinate system 
coincides with the linearization on every leaf of $W^1$. This implies that 
the linearizations $h^1_p$ depend smoothly on $p$.

Let $\psi_p$ be the restriction of $h_x$ to $W^1(p)$.
By  the construction of $h_x$,  the diffeomorphism $\psi_p$ can be expressed as
 $$
  \psi_p = \bar H_{0,\bar p}\circ h^1_x \circ H_{p,x},
 $$
where $H_{p,x}$ is the  holonomy map from $W^1 (p)$ to $W^1 (x)$ 
along the leaves of $W^2$, and $\bar H_{0,\bar p}$ is the projection
from $\Rn$ to $\Rn_{\bar p}=\Rn+{\bar p}$ along $\Rm$ in $\Rn\times \Rm$.
Since the holonomy map $H_{p,x}$ is affine by \fullref{smooth E},
 the map $\psi_p$ is also affine. 

We denote by $\bar h_p^1$ the coordinate representation of $h_p^1$, 
ie, 
 $$ 
  \bar h^1_p = d\psi_p(p) \circ h^1_p \circ \psi_p^{-1} \co
      \Rn_{\bar p} \cap U \to \Rn_{\bar p}.
 $$
Since $\psi_p$ is affine, we conclude that  $\bar h^1_p$ is also affine. 
We note that 
 $$
  \bar h^1_p(\bar p)= \bar p \; \text{ and }\;
  d\bar h^1_p(\bar p)=
  d\psi_p(p) \circ dh^1_p(p) \circ d(\psi_p^{-1}) (\bar p)=\Id,
  $$
since $dh^1_p(p)=\Id$ by \fullref{linearization}(2).
Hence the affine map  $\bar h_p^1$ is the identity map,
and therefore the maps $h^1_p$ depend uniformly $\Ci$ on $p$.
The uniformity comes from the fact that $h_x$ depends continuously
on $x$ in $\Ci$ topology.
\end{proof}


\subsection{Smooth affine connection} \label{connection}

In this section we construct an  $\a$--invariant $\Ci$ affine connection on $\M$.
We say that an affine connection $\nabla$   is of class $C^r$, 
$r\ge 0$, if $\nabla_X Y$ is a $C^r$ vector field for any two 
$C^\infty$ vector fields $X$ and $Y$. 

It was proved by Feres in \cite{F1} that an Anosov diffeomorphism which is 
1/2--pinched on its stable and unstable distributions preserves a unique 
continuous affine connection. In fact, such a connection exists on any 
invariant distribution where a diffeomorphism is a 1/2--pinched contraction
 \cite{F2}. The connections on the stable and unstable distributions can be
combined  into a connection on the manifold. However, the connection
on the (un)stable distribution is not known to be smooth transversally even 
if the distribution is smooth. Feres also noted in \cite{F2} that the exponential 
map of an invariant connection gives a nonstationary local linearization. 

We consider coarse Lyapunov distributions and use the nonstationary 
linearizations to obtain invariant affine connections on them. Their smoothness
will follow from the smoothness of the linearizations. It can be seen from \fullref{consistency} that these affine connections are the same as the ones defined
by the affine structures on the leaves of coarse Lyapunov foliations.

First we define an $\a$-invariant affine connection $\nabla^i$ along
each coarse Lyapunov foliation $W^i$.  At each point $x$ in $\M$, pull
back the affine connection $\bar \nabla^{E^i(x)}$ on the tangent space
$E^i (x)$ using the map $h^i_x\co W^i(x) \to E^i(x)$ to define an
$\a$--invariant affine connection $\nabla^i$ along each coarse
Lyapunov foliation $W^i$.  More precisely, for vector fields $X^i$ and
$Y^i$ on $W^i$ we define
 $$
  (\nabla^i_{X^i} Y^i)(x)=(h^i_x)_\ast^{-1} 
       \left( \bar \nabla^{E^i(x)}_{\bar X ^i} \,\bar  Y^i \right),
  $$ 
where $\bar X^i=(h^i_x)_\ast {X^i}$ and $ \bar Y^i =(h^i_x)_\ast Y^i$ are the 
push-forwards of $X^i$ and $Y^i$.
It is easy to see that $\nabla^i$ is an $\alpha$--invariant affine  connection,
which is as smooth as the dependence of $h^i_x$ on $x$. Thus, by
\fullref{smooth h}, $\nabla^i$ is uniformly $\Ci$ along the leaves of any
coarse Lyapunov foliation. It follows from \fullref{Ind Journe} that
 $\nabla^i$ is $\Ci$ on $\M$.

Now we define an affine connection $\nabla$ on $\M$ using a standard
construction. Let $X$ and $Y$ be two vector fields on $\M$. We decompose 
$X=\sum X^i$ and $Y=\sum Y^i$, where $X^i, Y^i \in E^i$. Then
 $$
   \nabla_X Y = \sum \nabla^i _{X^i} Y^i + \sum_{i \ne j} \Pi_j [X^i, Y^j],
 $$
where $\Pi_j$ is the projection onto $E^j$, defines an 
affine connection. Since the distributions $E^i$ and the connections $\nabla^i$
are $\Ci$, so is $\nabla$. Since $\nabla^i$ and $E^i$ are $\alpha$--invariant,
so is  $\nabla$. Thus we constructed an $\a$--invariant $\Ci$ affine connection on $\M$.


\subsection{Smooth conjugacy to a toral action} \label{smooth conjugacy}
 
In this section we complete the proof of \fullref{Main QC}. 
First we show that the action is $\Ci$--conjugate to an action 
by affine automorphisms of an infranilmanifold. 

Let $f$ be a transitive Anosov element of the action $\a$. 
Then $f$ preserves the  $\Ci$ affine connection $\nabla$ constructed in
\fullref{connection}. Also, its stable and unstable distributions are $C^\infty$ 
as direct sums of $\Ci$ coarse Lyapunov distributions.
It follows from the main result of  Y Benoist and F Labourie in \cite{BL} 
that $f$ is conjugate to an automorphism of an infranilmanifold $\N$ 
by a $\Ci$ diffeomorphism $\phi$.
It is known that any diffeomorphism commuting with an Anosov 
automorphism of an  infranilmanifold is an affine automorphism itself 
(see \fullref{Conjugacy}).  Hence $\phi$ conjugates $\a$  to an action 
$\bar \a$ by affine automorphisms of $\N$. 

Now we complete the proof by showing that $\N$ is finitely covered by a torus.
We use joint integrability of the coarse Lyapunov foliations and uniform 
quasiconformality of the action on them.

The conjugacy maps the (un)stable manifolds for elements of $\a$ to the 
(un)stable manifolds for elements of $\bar \a$, and the course Lyapunov 
foliations of $\a$ are mapped to the course Lyapunov foliations of $\bar \a$.
It follows that any two course Lyapunov foliations of $\bar \a$ are jointly
integrable.
Since the conjugacy is smooth, the Lyapunov functionals of the two actions
coincide, and $\bar \a$ is uniformly quasiconformal on its course Lyapunov 
distributions. Also, the action $\bar \a$ is TNS, and hence it has no
negatively proportional Lyapunov functionals.

The infranilmanifold $\N$ is finitely covered by a nilmanifold 
$N/\Gamma$, where $N$ is a simply connected nilpotent Lie group,
and $\Gamma$ is a cocompact lattice in $N$.
We need to show that $N$ is abelian. The Lie algebra $\n$ of $N$
splits into Lyapunov subspaces $\E_i$ with Lyapunov functionals $\chi_i$
of the action $\bar \a$. If for $u \in \E_i$ and $v \in \E_j$ the bracket $[u,v]$ is nonzero,
then $[u,v]$ belongs to a Lyapunov subspace with Lyapunov functional 
$\chi_i + \chi_j$. 
We recall that a coarse Lyapunov subspace is a direct sum of Lyapunov 
subspaces corresponding to positively proportional Lyapunov functionals. 

Suppose that $[u,v] \ne 0$ for some $u$ and $v$ which belong to two different 
coarse Lyapunov subspaces. Then $\chi_i$ and $\chi_j$ are not positively 
proportional and hence are not proportional. If follows that $\chi_i + \chi_j$ is 
not proportional to either one of
them. Hence $[u,v]$  belongs to a coarse Lyapunov subspace different from
the ones containing $u$ and $v$. This contradicts the fact that any two course 
Lyapunov foliations of $\bar \a$ are jointly integrable.

Suppose that  $u$ and $v$ are in the same coarse Lyapunov 
subspace $\E$. Uniform quasiconformality of the action on $\E$ implies that
all vectors in $\E$ are expanded/contracted at the same rate, ie,
there is only one Lyapunov functional $\chi$ corresponding to $\E$.
If $[u,v]$ is nonzero, then it belongs to a Lyapunov subspace with 
Lyapunov functional $2\chi$. But  this subspace must be  contained in $\E$,
which is impossible.

We conclude that $[u,v] = 0$ for any $u$ and $v$ in the Lie algebra $\n$. 
Thus the Lie group $N$ is abelian, and the infranilmanifold $\N$ is finitely 
covered by a torus. This completes the proof of \fullref{Main QC}.


\section[Proof of \ref{Main pinching} and its corollaries]{Proof of \fullref{Main pinching} and its corollaries} 
\label{Proof Main pinching}

\subsection[Proof of \ref{Main pinching}]{Proof of \fullref{Main pinching}} 

The proof of this theorem follows the same steps as the proof
of \fullref{Main QC}. Below we describe the necessary adjustments.

The coarse Lyapunov distributions are substituted by the distributions 
$E^i$ given in the theorem. The fact that any two such distributions are 
contracted by some element of the action is given by condition (2) of the
theorem. Since these distributions are intersections of
some stable distributions, they have most of the properties of coarse
Lyapunov distributions. In particular, their integral foliations are H\"older 
continuous with uniformly $\Ci$ leaves. 

Let $E^1$ and $E^2$ be two  distributions, and let $W^1$ and $W^2$ be 
the corresponding foliations. As in \fullref{regularity of sum}, we show 
that the leaves of the foliation $W=W^1 \oplus W^2$, given by assumption (1), 
are smooth. Then we show the smoothness of $E^1$ along the leaves of $W$ 
using the following proposition in place of \fullref{smooth distributions}. 

\begin{proposition} \label{smooth distributions pinching} 
Suppose that there exists an element $f$ that expands $E^1$ faster than it 
expands $E^2$. Then $E^1$ is $\Ci$  along the leaves of $W$, ie, its 
derivatives of all orders in $W$ directions exist, are continuous along the 
leaves of $W$ and uniformly bounded on $\M$.
\end{proposition}

\begin{proof}
We apply the $C^r$ Section Theorem as in the proof of
\fullref{smooth distributions}. The assumptions on $f$ directly imply 
that $\sup_{x \in M}  k_x <1$, ie, the induced map $F$ contracts the fibers. 
Moreover, since $f$ expands $W$, we see that $\a_x <1$, and thus 
 $$
   \sup_{x \in M} k_x \a_x ^r <1 \;\text{ for any } r>0.
 $$ 
The $C^r$ Section Theorem implies that the distribution $E^1$ is $C^\infty$ 
smooth along the leaves of $W$.
\end{proof}

As in \fullref{linearize}, we obtain nonstationary linearizations of the 
action along the leaves of  foliations $W^i$ and the associated 
affine structures. In this case, \eqref{(4.2)} follows from the 1/2--pinching 
assumption. 

It follows from condition (3) of the theorem and \fullref{smooth distributions pinching}  that the holonomies of $W^1$ 
between the leaves of $W^2$ are $\Ci$ along the leaves of $W=W^1 \oplus W^2$
and vise versa. Also, the 1/2--pinching assumption implies \eqref{(4.3)} with $\beta=1$. 
Hence we can apply \fullref{affine holonomies} to show that the holonomies 
are affine and uniformly $\Ci$ along 
the leaves of $W$. As in \fullref{holonomies}, 
it follows that all distributions $W^i$ are $\Ci$ on $\M$.

As in Sections \ref{smooth linearizations} and \ref{connection}, we show
that the nonstationary linearizations depend $\Ci$ on the base point and
construct a $\Ci$--invariant affine connection. As in the beginning of 
\fullref{smooth conjugacy}, we conclude that the action 
is $\Ci$--conjugate to an action by affine automorphisms of an infranilmanifold. 
\end{proof}


\subsection[Proof of Corollary \ref{TNS pinching}]{Proof of \fullref{TNS pinching}} 
\label{proof TNS pinching}

We need to show that  the infranilmanifold obtained in \fullref{Main pinching} is finitely covered by a torus. The argument is similar to 
the one in \fullref{smooth conjugacy}. We use
joint integrability of the coarse Lyapunov foliations and the 1/2--pinching. 
The only difference is in showing that  $[u,v] = 0$ for any $u$ and $v$ 
which belong to the same coarse Lyapunov subspace $\E$.  Suppose that 
 $[u,v] \ne 0$ for some $u$ and $v$ in $\E$.  Then the corresponding 
 Lyapunov functionals $\chi_i$ and $\chi_j$ are positively proportional.
 Hence they are  positively proportional to 
$\chi_i + \chi_j$, and the Lyapunov subspace of $\chi_i + \chi_j$ is 
also contained in $\E$. It is easy to see that this contradicts the existence 
of an element  which is 1/2--pinched on $\E$.


\subsection[Proof of \ref{functional pinching}]{Proof of \fullref{functional pinching}} 
\label{proof functional pinching}

It suffices to verify assumptions (2) and (3) of \fullref{Main pinching}. 
Let $E^1$ and $E^2$ be two coarse Lyapunov distributions. Let $\chi^1_{\min}$
and $\chi^2_{\min}$ be the corresponding functionals, and let $\H^1$ and $\H^2$
be their negative half-spaces in $\Rk$. If follows from assumption (2) 
of \fullref{functional pinching} that 
the negative half-space of $\chi^i_{\max}$ is $\H^i$, and that all elements in
$\H^i$ contract $E^i$.

Since the action is TNS Anosov, there exists an element which contracts both 
$E^1$ and $E^2$. Hence the intersection $\K=\H^1 \cap \H^2$ is a nonempty
open convex cone in $\Rk$.
It follows from condition (2) of the corollary that all elements in $\K$
contract both $E^1$ and $E^2$ and are 1/2--pinched on them. This verifies 
assumption (2) of \fullref{Main pinching}. 
To verify assumption (3) of the Theorem, we take an element 
$a$ in $\Zk \cap (-\K)$ close to the boundary of $\H^2$ and
away from the boundary of $\H^1$. Clearly, such an element expands $E^1$
faster than $E^2$. This completes the proof of the corollary. 


\subsection[Proof of \ref{perturbation}]{Proof of \fullref{perturbation}}

It is well known that an Anosov $\Zk$ action by toral automorphisms is
topologically conjugate to a $C^1$--small perturbation. To show the smoothness 
of the conjugacy we verify the assumptions  of \fullref{Main pinching} for
the perturbed action $\tilde \a$. 

First we consider the unperturbed action $\a$ by toral automorphisms.
For any coarse Lyapunov distribution $E$ we denote by $\chi_{\min}$ and 
$\chi_{\max}$ the minimal and the maximal of the Lyapunov functionals 
corresponding to $E$. The existence of a regular element which is 1/2--pinched
on $E$ implies that condition (2) of \fullref{functional pinching} is satisfied. 
Hence we can verify assumptions (2) and (3) of \fullref{Main pinching} 
for  $\a$ as in \fullref{proof functional pinching}. 

Now we consider a $C^1$--small perturbation $\ta$. For any coarse Lyapunov
distribution $E$ of $\a$ there exists a corresponding $\ta$--invariant distribution 
$\tE$ close to $E$. Each $\tE$ is the intersection of the corresponding stable
distributions for $\ta$, and the tangent bundle  is the direct sum of the
distributions $\tE$.
If the perturbation is sufficiently $C^1$--small, the assumptions (2) 
and (3) of \fullref{Main pinching} are satisfied for the distributions $\tE$.
Assumption (1) follows from the existence of topological conjugacy to $\a$.
Now \fullref{Main pinching} implies that the  conjugacy 
is $\Ci$.  
\endproof

\bibliographystyle{gtart}
\bibliography{link}

\end{document}